\documentclass[12pt]{article}

\usepackage[usenames,dvipsnames]{color}

\def\chi{{\mathcal{X}}}

\def\cald{{\mathcal{D}}}
\def\calx{{\mathcal{X}}}
\def\call{{\mathcal{L}}}

\def\({\left(}
\def\){\right)}

\def\vsp{\vspace*{1,5mm}\\ }

\def\bk{\bigskip }
\def\mk{\medskip }
\def\sk{\smallskip }

\def\n{\noindent }
\def\dd{\displaystyle}

\def\D{{\Delta}}

\def\barr{\begin{array}}
	\def\earr{\end{array}}

\def\bit{\begin{itemize}}

	\def\eit{\end{itemize}}

\def\D{{\Delta}}
\def\FP{Fokker--Planck}

\usepackage{amsmath, amsthm, amscd, amsfonts, amssymb}

\newtheorem{theorem}{Theorem}[section]
\newtheorem{proposition}[theorem]{Proposition}

\newtheorem{corollary}[theorem]{Corollary}

\newtheorem{lemma}[theorem]{Lemma}
\theoremstyle{definition}

\newtheorem{remark}[theorem]{Remark}

\def\1{^{-1}}
\def\vsp{\vspace*{2mm}\\ }

\def\calf{{\mathcal{F}}}

\def\calx{{\mathcal{X}}}

\def\rr{{\mathbb{R}}}
\def\rrd{{\mathbb{R}^d}}
\def\nn{{\mathbb{N}}}
\def\9{{\infty}}

\def\lbb{{\lambda}}

\def\b{{\beta}}

\def\ov{\overline}
\def\vf{{\varphi}}

\def\ooo{{\Omega}}
\def\pp{{\partial}}
\def\D{{\Delta}}
\def\vp{{\varepsilon}}
\def\barr{\begin{array}}
	\def\earr{\end{array}}

\def\dd{\displaystyle}
\def\bk{\bigskip }
\def\sk{\smallskip}
\def\n{\noindent }

\def\vsp{\vspace*{2mm}\\ }
\def\ff{\forall }
\def\({\left(}
\def\){\right)}
\def\<{\left<}
\def\>{\right>}

\title{Nonlinear Fokker--Planck    equations with~time-dependent coefficients}
\author{Viorel Barbu\thanks{Octav Mayer Institute of Mathematics of  Romanian Academy,     Ia\c si, Romania.  Email: vbarbu41@gmail.com}\and Michael R\"ockner\thanks{Fakult\"at f\"ur Mathematik, Universit\"at Bielefeld,  D-33501 Bielefeld, Germany.  Email: roeckner@math.uni-bielefeld.de}}
\date{}

\begin{document}
	\maketitle
	\begin{abstract}
		\n An   operatorial based approach is used here to prove the existence and uniqueness of a strong  solution $u$ to  the time-varying nonlinear  Fokker--Planck equation\vspace*{-2mm}
		$$\barr{l}
		u_t(t,x)-\Delta(a(t,x,u(t,x))u(t,x))+{\rm div}(b(t,x,u(t,x))u(t,x))=0\\\hfill\mbox{ in }(0,\9)\times\rr^d,\vspace*{-2mm}\\ u(0,x)=u_0(x),\ x\in\rr^d\earr$$in the Sobolev space $H\1(\rr^d)$, under appropriate conditions on the  $a:[0,T]\times\rr^d\times\rr\to\rr$ and $b:[0,T]\times\rr^d\times\rr\to\rr^d.$  Is is proved also that, if  $u_0$ is a density of a probability measure, so is $u(t,\cdot)$  for all \mbox{$t\ge0$}. Moreover, we construct a weak   solution to the McKean-Vlasov SDE associated with the  Fokker-Planck equation such that $u(t)$ is the density of its time marginal law.\sk
	\end{abstract}
	\section{The problem}\label{s1}
	We  shall study here the existence and uniqueness of a strong solution in the~Sobolev space $H\1(\rr^d)$ to the nonlinear time-varying \FP\ equation\newpage	
	\begin{equation}\label{e1.1}
	\barr{l}
	u_t(t,x)-\Delta(a(t,x,u(t,x))u(t,x))+{\rm div}(b(t,x,u(t,x))u(t,x))=0\\\hfill\mbox{ on }(0,T)\times\rr^d,\\
	u(0,\cdot)=u_0,
	\earr\end{equation}where $d\ge1$ and $0<T<\9$.
	In the following, we   use the notations
	\begin{equation}
	\beta(t,x,r)\equiv a(t,x,r)r,\ b^*(t,x,r)\equiv b(t,x,r)r,\ (t,x,r)\in[0,T]\times\rr^d\times\rr,
	\end{equation}and we shall assume that
	\begin{itemize}
		\item[\rm(H1)] $\beta\in C^1([0,T]\times\rr^d\times\rr)$, $\beta,\beta_r$ are bounded and, for all $t,s\in[0,T],$ $  x\in\rr^d,\ r,\bar r\in\rr,$
		\begin{equation}
		(\beta(t,x,r)-\beta(t,x,\bar r))(r-\bar r)\ge\nu|r-\bar r|^2,\label{e1.3}\end{equation}
		\begin{equation}\label{e1.3prim}
		|\beta_r(t,x,r)-\beta_r(s,x,r)|\le h(x)|t-s|\beta_r(t,x,r),
		\end{equation}
		\begin{equation}
		\barr{r}
		 |\beta(t,x,r)-\beta(s,x,r)|+|\beta_x(t,x,r)-\beta_x(s,x,r)|\\\le h(x)|t-s|(1+|r|),\earr\label{e1.4}
		\end{equation}
		\begin{equation}
		|\beta_t(t,x,r)|+|\beta_x(t,x,r)|\le h(x)|r|,\  \ff(t,x,r)\in[0,T]{\times}\rr^d{\times}\rr.\label{e1.4az}
		\end{equation}
		\item[\rm(H2)] $r\to b(t,x,r)$ is $C^1$,  $\ff(t,x)\in[0,T]\times\rr^d$, $b,rb_r$ are bounded,  and $\ff t,s\in[0,T],$ $ x\in\rr^d,$ $ r\in\rr,$
		\begin{equation}\label{e1.4a}
		|b^*(t,x,r)-b^*(s,x,r)|\le h(x)|t-s|(1+|b^*(t,x,r)|),\end{equation}
		and
		\begin{equation}
		\label{e1.7prim}
		|b^*(t,x,r)|\le h(x)|r|,
		\end{equation}
	\end{itemize}
	where $\nu>0$, $h\in L^\9(\rr^d)\cap L^2(\rr^d)$, $h\ge0$. Here   $b_r=\frac{\pp b}{\pp r}$, \mbox{$\beta_x=\nabla_x\beta$,}  $\beta_r=\frac\pp{\pp r}\,\beta$ and div, $\D, \nabla$ are taken with respect to the spatial variable \mbox{$x=(x_j)^d_{j=1}$} and are considered in the sense of Schwartz distributions on $\rr^d$. Here, we note that by \eqref{e1.3} we have $\beta_r\ge\nu>0.$
	
	Equation \eqref{e1.1} is relevant in the description   of particle transport in disor\-dered media (so called anomalous diffusions) with time dependent temperature. For instance, if $\beta(t,x,r)\equiv\gamma(t,x)\ln(1+|r|)$, \eqref{e1.1} describes the classical bosons dynamics. (Other significant examples are found in \cite{7}--\cite{6a}.)   In appli\-cations, the solution $u$ to \eqref{e1.1} is a probability density and so the main  interest is to find distributional solutions $u$ which satisfy the  conditions 
		
	\begin{eqnarray}
	&u(t,x)\ge0,\ \ff t\in[0,T] \mbox{ and a.e. }x\in\rr^d,\label{e1.5}\\
	&\dd\lim_{t\to t_0}\int_{\rr^d}\!\!u(t,x)\psi(x)dx=\int_{\rr^d}\!\!u(t_0,x)\psi(x)dx,\ \ff t_0\in[0,T],\ \psi\in C_b(\rr^d),\qquad \label{e1.6}\\[1mm]
	&\dd\int_{\rr^d}u(t,x)dx=1,\, \ff t\in[0,T].\label{e1.7}
	\end{eqnarray}
	This means that $u=u(t,x)$ is a probability density (in $x$) and is weak-star continuous in $t\in[0,T]$ in the space of bounded measures in  $(C_b(\rr^d))^*$ on $\rr^d$. In equation \eqref{e1.1}, $b=b(t,x,r)$ is the drift coefficient which represents the field force acting on the particle, while $a=a(t,x,r)$ is the diffusion coefficient.  Indeed, if $u$ is a solution to \eqref{e1.1} satis\-fying \eqref{e1.5}--\eqref{e1.7}, then one can find a pro\-ba\-bi\-listic representation  of $u$ as the time marginal law of the solution $X=X(t)$, $t\in[0,T]$, to the  McKean-Vlasov stochastic differential equation (of Nemytskii-type)
	\begin{equation}\label{e1.8}
	\hspace*{-11mm}\barr{l}
	\dd dX(t){=}b\!\(\!t,X(t),\frac{d\call_{X(t)}}{dx}\,(X(t))\!\)\!dt{+}\sigma\!\(\!t,X(t),\frac{d\call_{X(t)}}{dx}\,(X(t))\!\)\!dW(t),\\
	X(0){=}\xi_0,\earr\hspace*{-11mm}\end{equation}where $\sigma=\sqrt{2a},$ $\call_{X(t)}$ is the law of $X(t)$, $W=\{W_j\}^d_{j=1}$ is an $(\calf_t)$-Brownian motion on a filtered probability space $(\ooo,\calf,(\calf_t),\mathbb{P})$ and $\xi_0:\ooo\to\rr^d$ is\break $\calf_0$-measurable such that $\mathbb{P}\circ\xi^{-1}_0(dx)=u_0(x)dx.$ In fact, using the superposition principle for linear \FP\ equations (see \cite{6}, which in turn generalizes \cite{12a}), it was shown in \cite{3a} %(based on \cite{12a} and \cite{6}), 
	that there is a {\it probabilistically weak solution} to \eqref{e1.8} such that
	\begin{equation}\label{e1.9}
	u(t,x)dx=\mathbb{P}\circ(X(t))^{-1}(dx),\ \ff t\ge0.\end{equation}(See Corollary \ref{c2.2} below.) Thus, the results of this paper and \cite[Section 2]{3a} lead to the first proof of existence of weak solutions to the McKean--Vlasov SDE \eqref{e1.8}.

The existence of a weak solution $u$ for \eqref{e1.1}  was proved in \cite{2}, \cite{3}, \cite{4a}  and, for a more general class of \FP\ equations in \cite{3a}, in  the autonomous case $\beta(t,r)\equiv\beta(r)$, $b(t,x,u)\equiv b(x,u)$. In this case,   the generation theorem for nonlinear semigroups of contractions in $L^1(\rr^d)$  was the essential instrument to prove the well posedness of \eqref{e1.1}.  As it is well known, the nonlinear semigroup theory does not extend completely to the Cauchy problem for differential equations with time-dependent accretive operators, and so here we shall use a different approach based on the general existence theory of the time-dependent Cauchy problem in Banach spaces (\cite{5}) combined with an $L^1(\rr^d)$ analysis of nonlinear \FP\ equations similar to that de\-ve\-loped in the above mentioned works \cite{3}, \cite{3a}, \cite{4a}.
	
	The literature on both (nonlinear) \FP\ equations and McKean--Vlasov SDEs is huge. We here refer to \cite{9prim} and \cite{9secund}, respectively, as well as the references therein. We would like to point out, however, that in the literature the coefficients are usually assumed to be weakly continuous in its measure, i.e., time marginal law, dependence, whereas in \eqref{e1.8} this is not the case. The measure dependence is via the time marginal law density evaluated at $t$ and $x$, $X(t)$ respectively ({\it Nemytskii-type dependence}).

 \bk
 \noindent{\bf Notation.} For $1\le p\le\9$, denote by $L^p(\rr^d)$  (or $L^p$) the space of $p$-Lebesgue summable functions on $\rr^d$, with norm denoted $|\cdot|_p$. By $L^p_{\rm loc}(\rr^d)=L^p_{\rm loc}$ we denote the corresponding local space on $\rr^d$. The scalar \mbox{product} in $L^2(\rr^d)$ is denoted by $\<\cdot,\cdot\>_2$. By $\cald'(\rr^d)$ and $\cald'((0,T)\times\rr^d)$ we denote the space of Schwartz distri\-butions on $\rr^d$ and $(0,T)\times\rr^d,$ respectively. By $H^k(\rr^d)$, de\-no\-ted $H^k$, $k{=}1,2,$ the Sobolev space $\{u\in L^2(\rr^d);\, D^\alpha u\in L^2(\rr^d),\, |\alpha|\le k\}$, where $\alpha=(\alpha_1,...,\alpha_d),$ $|\alpha|=\alpha_1+\alpha_2+\cdots+\alpha_d$, $D^\alpha =\frac{\pp^{\alpha_1+\alpha_2+\alpha_d}}{\pp x^{\alpha_1}_1...\pp x^{\alpha_d}_d}$
	and the derivatives are taken in the sense of $\cald'(\rr^d)$. By $H^{-1}(\rr^d)$, also written $H\1$, we denote the dual space of $H^1(\rr^d)$ with the scalar product
	\begin{equation}\label{e1.10}
	\<u,v\>_{-1,\vp}=\int_{\rr^d}(\vp I-\Delta)^{-1}uv\,dx,\ \ff u,v\in H^{-1}(\rr^d),\end{equation}
	and the corresponding norm $|\cdot|_{-1,\vp}.$ Here $(\vp I-\Delta)^{-1}u=y\in H^2(\rr^d)$ is~defined for $\vp>0$ and $u\in L^2(\rr^d)$ by the equation
	\begin{equation}\label{e1.11}
	\vp y-\Delta y=u\mbox{ in }\cald'(\rr^d).\end{equation}
	$H^{-2}(\rr^d)$ is the dual space of $H^2(\rr^d)$ with the norm $|u|_{-2}=|(\vp I-\Delta)^{-1}u|_2,$ $u\in H^{-2}(\rr^d).$ By $W^{1,p}([0,T];H^{-1}(\rr^d))$, $1\le p\le\9$ we denote the\break space of absolutely continuous functions $u:[0,T]\to H^{-1}(\rr^d)$ such that\break $\frac{du}{dt}\in L^p(0,T;H^{-1}(\rr^d))$.
	Here $\frac{du}{dt}$ is the derivative of $u$ in the sense of $H\1(\rr^d)$-valued distributions on $(0,T)$ or, equivalently, a.e. on $(0,T)$.
	
	Denote by $C([0,T];L^2(\rr^d))$ the space of continuous $L^2(\rr^d)$-valued functions and by $C_w([0,T];L^2(\rr^d))$ the space of weakly continuous functions\break $u:[0,T]\to L^2(\rr^d)$. Finally, $C_b(\rr^d)$ and $C_b([0,T]\times\rr^d)$ are the space of all  continuous and bounded functions on $\rr^d$ and $[0,T]\times\rr^d$, respectively.
	
	The main existence result for equation \eqref{e1.1} will be formulated in Section~2 and proved in Section 3.
	%The periodic problem \eqref{e1.11a}-\eqref{e1.11aa} will be treated in Section 4.
	
	\section{The main existence result}\label{s2}
	\setcounter{equation}{0}
	
	\begin{theorem}\label{t2.1} Assume that Hypotheses {\rm(H1), (H2)} hold. Let $0<T<\9$ be arbitrary, but fixed. Let $D_0:=\{u\in L^2(\rr^d);\,\beta(0,\cdot,u)\in H^1(\rr^d)\}$. Then, for each $u_0\in D_0$,   there is a unique strong solution $u$ to \eqref{e1.1} satisfying 	
		\begin{eqnarray}
		&u\in C([0,T]; L^2(\rr^d))\cap W^{1,2}([0,T];H^{-1}(\rr^d)),\label{e2.1}\\
			& \label{e2.2}
		u,\beta(\cdot,u)\in L^2(0,T;H^1(\rr^d)).	\end{eqnarray}
		If $u_0\in L^1(\rr^d)\cap L^2(\rr^d),$ then
		$u\in L^\9(0,T;L^1(\rr^d))$. Moreover, it satisfies the weak continuity condition \eqref{e1.6} and $|u(t)|_1\le |u_0|_1$, for all $ t\in [0,T].$
		
		In addition, if $u_0\in L^1(\rr^d)\cap D_0$ is 
		such that\vspace*{-2mm}
		\begin{equation}\label{e2.3}
		u_0\ge0,\ \mbox{ or }\ \  \int_{\rr^d}u_0(x)dx=1,\vspace*{-2mm}\end{equation}then $u$ satisfies also conditions \eqref{e1.5}, \eqref{e1.7}, respectively.
		
		Finally, assume  that \eqref{e1.7prim} is replaced by the stronger condition
		\begin{equation}
		\label{e2.3prim}
		|b^*(t,x,r)-b^*(t,x,\bar r)|\le h(x)|r-\bar r|,\ \ff\,t\in[0,T],\ x\in\rr^d,\ r,\bar r\in\rr.
		\end{equation} Then, for all $u_0,\bar u_0\in L^1(\rr^d)\cap D_0$, the corresponding solutions $u(t,u_0),$ $u(t,\bar u_0)$ to \eqref{e1.1} satisfy
		\begin{equation}\label{e2.3a}
		|u(t,u_0){-}u(t,\bar u_0)|_1\le|u_0{-}\bar u_0|_1,\ \ff t\in[0,T].\end{equation}Moreover, if  $u_0\in L^1(\rr^d)\cap D_0\cap L^\9$ and \begin{equation}\label{e2.5prim}
		\Lambda(b,\beta)=\sup\{|b_x(t,x,r)r|{+}|\Delta\beta(t,x,r)|; t\!\in\![0,T], x\!\in\!\rr^d, r\!\in\!\rr\}<\9,\end{equation} then $u\in L^\9((0,T)\times\rr^d)$, and $\|u\|_\9\le\Lambda(b,\beta)T+\|u\|_\9$.
	\end{theorem}
	Here, as above, $\beta(t,x,r)\equiv a(t,x,r)r$.

	By a strong solution to equation  \eqref{e1.1} we mean a function $u$ satisfying \eqref{e2.1}, \eqref{e2.2} and equation \eqref{e1.1} in $H^{-1}(\rr^d)$, a.e.  in $t\in(0,T)$ with $\Delta$ and div taken in the sense of Schwartz distributions on $\rr^d$. (We note that, by \eqref{e2.1}, \eqref{e2.2}, we have
	$u_t\in L^1(0,T;H^{-1}(\rr^d)),$ $\Delta\beta(\cdot,u)\in L^2(0,T;H^{-1}(\rr^d)),$ ${\rm div}(b(\cdot,u)u)\in L^2(0,T;H^{-1}(\rr^d))$.) In other words, $u\in L^2(0,T;L^2(\rr^d))$ is a solution to \eqref{e1.1} in the sense of $H\1(\rr^d)$-valued vectorial distributions on~$(0,T)$.
	
	It is easily seen that the solution $u$ given by Theorem \ref{t2.1} is also a solution in the sense of Schwartz distributions on $(0,T)\times\rr^d$, that is,
		
	\begin{equation}\label{e2.4}
	\barr{l}
	 \dd\int^T_0\int_{\rr^d}(u(t,x)\cdot\vf_t(t,x)+a(t,x,u(t,x))u(t,x)\Delta\vf(t,x)\\
	+u(t,x)b(t,x,u(t,x))\cdot\nabla\vf(t,x))dt\,dx
	+\dd\int_{\rr^d}\vf(0,x)u_0(x)dx=0,\vsp
	\hfill \ff\vf\in C^2_0([0,T)\times\rr^d),\earr\end{equation}where $C^2_0([0,T)\times\rr^d)$ is the space of all twice  differentiable functions with compact support in $[0,T)\times\rr^d$.
	In particular, by \eqref{e2.3a} the  $L^1$-stability   of the solution $u=u(t,u_0)$ follows and also that $t\to u(t,u_0)$ is a contraction flow in the space $L^1(\rr^d)$. 
	
 \begin{remark}\label{r2.1prim} \rm
 	By \eqref{e2.3a}, it follows  via a density argument that, in fact, for every {$u_0\in L^1$,} there is a   distributional solution $u$ in the sense of \eqref{e2.4} satisfying \eqref{e2.3a}. It follows by \cite[Lemma 8.1.2]{0} that $t\mapsto u(t,x)dx$ has a weakly continuous $dt$-version on $[0,T]$, which thus  satisfies \eqref{e1.6}.\end{remark}
As noted earlier,   by virtue of Theorem \ref{t2.1}, we have the following consequence (see Section 3 below for its proof):
 \begin{corollary}\label{c2.2} Assume that {\rm(H1), (H2)} and $u_0\in L^1\cap D_0$ such that \eqref{e2.3} holds. Then, there exists a $($probabilistically$)$ weak solution $X$ to the McKean--Vlasov SDE \eqref{e1.8}, with $\sigma=\sqrt{2a}$, where $a$ is as in {\rm(H1)}, such that  the strong solution $u$ to \eqref{e1.1} given by Theorem {\rm\ref{t2.1}} is its time marginal law density. In particular, we have the representation formula \eqref{e1.9} for $u$. The same is true for all $u_0\in L^1$ if   \eqref{e2.3prim} holds.
 \end{corollary}	
	We shall prove Theorem \ref{t2.1} in several steps and the first one is an existence result for an approximating equation corresponding to \eqref{e1.1}.	 
	\begin{proposition}\label{p2.4} Assume that {\rm(H1), (H2)} hold and that
$ u_0\in D_0.$ Then, for each $\vp>0$, the equation
 \begin{equation}\label{e2.6}
 \hspace*{-4mm}\barr{l}
 u_t-\Delta\beta(t,x,u)+\vp\beta(t,x,u)+{\rm div}(b(t,x,u)u)=0,\, t\in(0,T),\, x\in\rr^d,\vsp
 u(0,x)=u_0(x),\earr\!\!\!\!\end{equation}has a unique strong solution $u=u_\vp\in W^{1,\9}([0,T];H\1(\rr^d))\cap L^2(0,T;H^1(\rr^d))$, such that $\beta(u_\vp)\in L^2(0,T;H^1(\rr^d))$. In particular, $u_\vp\in C([0,T];L^2(\rr^d))$.
\end{proposition}
Then the first part of Theorem \ref{t2.1} will follow from Proposition \ref{p2.4} by letting $\vp\to0$ in equation \eqref{e2.6}.

As regards Proposition \ref{p2.4}, it will be derived from an abstract existence result for the Cauchy problem in the space $H\1(\rr^d)$
 \begin{equation}\label{e2.7}
 \barr{l}
 \dd\frac{du}{dt}+A_\vp(t)u=0,\ t\in(0,T),\vsp
 u(0)=u_0,\earr\end{equation}where $A_\vp(t):D(A_\vp(t))\subset H^{-1}(\rr^d)\to H^{-1}(\rr^d)$ is for each $t\in[0,T]$ a quasi-$m$-accretive operator in $H^{-1}(\rr^d)$ to be defined below.

 \section{Proofs}\label{s3}
 \setcounter{equation}{0}
 We shall prove first Proposition \ref{p2.4}. To this purpose, we define for $t\in[0,T]$ and $\vp>0$ in the space $H^{-1}(\rr^d)$ the operator $A_\vp(t):D(A_\vp(t))\subset H^{-1}(\rr^d)\to H^{-1}(\rr^d)$
 \begin{equation}\label{e3.1}
 \barr{rl}
 (A_\vp(t)y)(x)\!\!\!&=-\Delta\beta(t,x,y(x))+\vp\beta(t,x,y(x))\vsp
 &+\,{\rm div}(b(t,x,y(x))y(x)),\ \ff y\in D(A(t)),\ x\in\rr^d,\vsp
 D(A_\vp(t))&=\{y\in L^2(\rr^d);\beta(t,\cdot,y)\in H^1(\rr^d)\},\ t\in[0,T],\earr\end{equation}where the differential operators $\Delta $ and div are taken in the sense of the Schwartz distribution space  $\cald'(\rr^d)$.

  By (H1) it follows immediately that $$H^1(\rr^d)\cap L^\9(\rr^d)\subset D(A_\vp(t))\subset H^1(\rr^d),$$ hence $D(A_\vp(t))$ is dense in $L^2(\rr^d)$.

  By Hypotheses (H1), (H2), it is easily seen that, setting $D_t:=D(A_\vp(t))$,
  $$A_\vp(t)(D_t)\subset H\1(\rr^d), \ \ff t\in[0,T].$$
 Since below we fix $\vp>0$, for simplicity we write $\<\cdot,\cdot\>_{-1}$ and $|\cdot|_{-1}$ instead of
 $\<\cdot,\cdot\>_{-1,\vp}$, $|\cdot|_{-1,\vp}$, respectively. We also have
 \begin{lemma}\label{l3.1} Let $t\in[0,T]$ and $\lbb_0:=\frac{2\nu}{|b^*_r|^2_\9}.$ Then, the operator $A_\vp(t)$ is quasi-$m$-accretive in $H\1(\rr^d)$, that is, $(I+\lbb A_\vp(t))\1:H\1(\rr^d)\to H\1(\rr^d)$ is single-valued and 
 		
 \begin{equation}\label{e3.2}
 \|(I+\lbb A_\vp(t))^{-1}\|_{{\rm Lip}_\vp(H^{-1}(\rr^d))}\le\(1-\frac\lbb{\lbb_0}\)^{-1},\ \ff\lbb\in(0,\lbb_0), \end{equation}where $\|\cdot\|_{{\rm Lip}_\vp}$ means Lipschitz norm w.r.t. the norm $|\cdot|_{-1,\vp}$ on $H\1.$ 

 In~particular, $(A_\vp(t),D_t)$ is closed for every $t\in[0,T]$.

Moreover, $D_t=D_0$, for every $t\in[0,T]$, and  for all $\lbb\in(0,\lbb_0)$ we~have
 \begin{equation}\label{e3.3}
\barr{l}
 |(I+\lbb A_\vp(t))^{-1}u-(I+\lbb A_\vp(s))^{-1}u|_{-1,\vp}\vsp\qquad\le \lbb|t-s|L(|u|_{-1,\vp})(1+|A_\vp(t)u|_{-1,\vp}),\  \ff u\in D_0,\ s,t\in[0,T],\earr\end{equation}where $L:[0,\9)\to[0,\9)$ is a monotone nondecreasing function.
\end{lemma}

\noindent{\bf Proof.} One must prove that there exists $\lbb_0>0$ such that, for each $t$ and $\lbb\in (0,\lbb_0)$, the operator $(I+\lbb A_\vp(t))\1$ is single-valued (on its domain) and that \eqref{e3.3} holds. 
In the following, we simply write
$$|\cdot|_{-1,\vp}=|\cdot|_{-1},\ \<\cdot,\cdot\>_{-1,\vp}=\<\cdot,\cdot\>_{-1}.$$

Let $v$ be arbitrary in $H^{-1}=H\1(\rr^d)$ and consider the equation
\begin{equation}
\label{e3.3prim}
u+\lbb A_\vp(t)u=v.
\end{equation}
In the following, we shall simply write $$\mbox{$\beta(t,x,u)=\beta(t,u),$ $b(t,x,u)=b(t,u)$, $b^*(t,x,u)=b^*(t,u)$.}$$
Let us consider \eqref{e3.3prim} in $\cald'(\rr^d)$, i.e., we seek a function $u\in L^2$ such that
 \begin{equation}\label{e3.4}
 u-\lbb(\Delta-\vp I)(\beta(t,u))+\lbb {\rm div}(b^*(t,u))=v\mbox{ in }\cald'(\rr^d).\end{equation}Equivalently,
 \begin{equation}\label{e3.5}
 (\vp I-\Delta)^{-1}u+\lbb\beta(t,u)+\lbb(\vp I-\Delta)^{-1}{\rm div}(b^*(t,u))=(\vp I-\Delta)^{-1}v.\end{equation}
 Let $u,\bar u\in L^2$. Then
 $$|u-\bar u|^2_{-1}=\vp|(\vp I-\Delta)\1(u-\bar u)|^2_2+|\nabla(\vp I-\Delta)\1(u-\bar u)|^2_2.$$Hence, by (H2) we have
\begin{equation}\label{e3.5prim}
\barr{l}
\<(\vp I-\Delta)^{-1}{\rm div}(b^*(t,u)-b^*(t,\bar u)),u-\bar u\>_2\vsp
\quad\qquad\ge
 -|b^*_r|_\9|\nabla(\vp I-\D)\1(u-\bar u)|_{2}|u-\bar u|_2\vsp
\quad\qquad\ge-|b^*_r|_\9|u-\bar u|_{-1}|u-\bar u|_2 \ge-\dd\frac{|b^*_r|^2_\9}{2\nu}\,|u-\bar u|^2_{-1}-\frac\nu2\,|u-\bar u|^2_2,
\earr\end{equation}
and, by \eqref{e1.3},
$$\<\beta(t,u)-\beta(t,\bar u),u-\bar u\>_2\ge\nu|u-\bar u|^2_2.$$
This  implies that the operator $F:L^2\to L^2$,
$$Fu=(\vp I-\Delta)^{-1}u+\lbb\beta(t,u)+\lbb(\vp I-\Delta)^{-1}{\rm div}(b^*(t,u)),\ u\in L^2,$$satisfies the condition
  \begin{equation}\label{e3.6}
 \dd\<Fu-F\bar u,u-\bar u\>_2\ge\frac{\nu\lbb}2\,|u-\bar u|^2_2+\(1-\frac{\lbb|b^*_r|^2_\9}{2\nu}\)|u-\bar u|^2_{-1},\vsp \ff u,\bar u\in L^2.\end{equation}Hence, for $\lbb\in(0,\lbb_0)$,  the operator $F:L^2\to L^2$ is accretive, continuous and coercive. Therefore, it is surjective (see, e.g., \cite{1}, p.~104). So, \eqref{e3.5} has a solution $u\in L^2$ and so, by \eqref{e3.4},  $\beta(t,u)\in H^1$. So, $u$ solves \eqref{e3.3prim}.  By \eqref{e3.6}, we also have
  $$|F^{-1}(v)-F^{-1}(\bar v)|_2\le\frac2{\nu\lbb}\,|v-\bar v|_2, \ \ff v,\bar v\in L^2,$$for $0<\lbb\le\lbb_0.$

 Moreover, by \eqref{e3.5}--\eqref{e3.6}, we see that, for $v,\bar v\in H^{-1}$ and the corres\-pon\-ding solutions $u,\bar u\in L^2$ for \eqref{e3.6}, we have
 $$\frac{\nu\lbb}2\,|u-\bar u|^2_2+\(1-\frac{\lbb}{\lbb_0}\)|u-\bar u|^2_{-1}
 \le|u-\bar u|_{-1}|v-\bar v|_{-1}.$$Therefore, $(I+\lbb A_\vp(t))\1$ is single valued and
 $$ |(I+\lbb A_\vp(t))^{-1}v-(I+\lbb A_\vp(t))^{-1}\bar v|_{-1}\le\(1-\dd\frac{\lbb}{\lbb_0}\)^{-1}|v-\bar v|_{-1},\, \ff v,\bar v\in  H^{-1},$$for $0<\lbb<\lbb_0.$ This implies \eqref{e3.2}, as claimed.

 Let us show now that $D(A_\vp(t))$ is independent of $t$. Indeed, by \eqref{e1.3prim} we~have
\begin{equation}\label{e3.8prim}
 \beta_r(s,x,r)\le(h(x)|t-s|+1)\beta_r(t,x,r),\,
 \ff \,t,s\in[0,T],\, x\in\rr^d,\, r\in\rr.\end{equation}Let $y\in L^2$. Then, by \eqref{e1.4}, $y\in D(A_\vp(t))$ if and only if   $\beta_r(t,\cdot,y)\nabla y\in L^2$. Hence, by \eqref{e3.8prim} we also have
 $$\beta_r(s,y)\nabla y\in L^2,\ \ \ff\,s\in[0,T],$$that is, $y\in D(A_\vp(s))$, $\ff \,s\in[0,T]$, as claimed.

 Therefore, $D(A_\vp(t))$ is independent of $t$, and
 $$D(A_\vp(t))=D_0=\{u_0\in L^2(\rr^d);\beta(0,u_0)\in H^1\}.$$

  To prove \eqref{e3.3}, we note that, by Lemma 3.2 in \cite{5}, \eqref{e3.3} holds if
  {\begin{equation}\label{e3.9a}
  \barr{r}
  |A_\vp(t)u-A_\vp(s)u|_{-1}\le C|t-s|L(|u|_{-1})(1+|A_\vp(t)u|_{-1}),\vsp \ff u\in D_0,\ s,t\in[0,T],\earr\end{equation}where $L:[0,\9)\to[0,\9)$ is a nondecreasing monotone  function.
 To this purpose, we note that by \eqref{e3.1} and \eqref{e1.10} we have  
  \begin{equation}
  \label{e3.9b}
  \barr{l}
  |A_\vp(t)u-A_\vp(s)u|^2_{-1}=
 \left<\beta(t,u)-\beta(s,u)\right.\vsp
  \qquad+(\vp I-\Delta)\1{\rm div}(b^*(t,u)-b^*(s,u)),
\vp(\beta(t,u)-\beta(s,u))\vsp
\qquad
  \left.-\Delta(\beta(t,u)-\beta(s,u))+{\rm div}(b^*(t,u)-b^*(s,u))\right>_2\vsp
  \qquad
  =|\nabla(\beta(t,u)-\beta(s,u))|^2_2
  +|{\rm div}(b^*(t,u)-b^*(s,u))|^2_{-1}\vsp
  \qquad
  +\vp|\beta(t,u)-\beta(s,u)|^2_2
  \vsp\qquad
  +2\left<\beta(t,u)-\beta(s,u),{\rm div}(b^*(t,u)-b^*(s,u))\right>_2,\vsp
  \qquad
  \hfill \ff\,s,t\in[0,T],\ u\in D_0.
  \earr
  \end{equation}
  We also have
$$|{\rm div}\,f|_{-1}\le C|f|_2,\ \ \ff\,f\in  L^2,$$and so, by \eqref{e1.4a}, we have
\begin{equation}\label{e3.11a}
	\hspace*{-3mm}\barr{ll}
 |{\rm div}(b^*(t,u)-b^*(s,u))|_{-1}\!\!\!&\le C|b^*(t,u)-b^*(s,u)|_2\vsp&
 \le C|t-s|(1+|b^*(t,u)|_2),\  \ff(s,t)\in(0,T),\earr\hspace*{-3mm}\end{equation}By \eqref{e1.3}--\eqref{e1.4az} it follows that
$$\barr{l}
|\beta(t,u)-\beta(s,u)|+|\nabla(\beta(t,u)-\beta(s,u))|\vsp 
\qquad\le|(\beta_r(t,u)-\beta_r(s,u))\nabla u|+|\beta_x(t,u)-\beta_x(s,u)|+|\beta(t,u)-\beta(s,u)|\vsp 
\qquad\le2h|t-s|(|\nabla\beta(t,u)|+
|\beta_x(t,u)|+|u|+1),\mbox{ a.e. in }(0,T)\times\rrd.
\earr$$This yields
\begin{equation}\label{e3.11aa}
	\barr{l}
 |\nabla(\beta(t,u)-\beta(s,u))|_2+|\beta(t,u)-\beta(s,u)|_2\vsp\qquad\le C|t-s|(1+|\beta(t,u)|_2+|\nabla\beta(t,u)|_2),\  \ff s,t \in[0,T].\earr\end{equation}Then, by \eqref{e3.9b}--\eqref{e3.11aa}, it follows that   
  \begin{equation}
  \label{e3.9bb}
  \barr{l}
  |A_\vp(t)u-A_\vp(s)u|_{-1}
  \le C|t-s|(1+|\beta(t,u)|_2\vsp
  \qquad+|\nabla\beta(t,u)|_2+|b^*(t,u)|_2),\  \ff\,t,s\in[0,T],\ u\in D_0.\earr
  \end{equation}
  On the other hand, we have as in \eqref{e3.9b} and using \eqref{e1.4az},  
   \begin{equation}
  \label{e3.9bbb}
  \hspace*{-5mm}\barr{ll}
  |A_\vp(t)u|^2_{-1}\!\!&
  =|\nabla\beta(t,u)|^2_2+|{\rm div}\,b^*(t,u)|^2_{-1}\vsp
  & +\,2\left<\beta(t,u),{\rm div}\,b^*(t,u)\right>_2+\vp|\beta(t,u)|^2_2\vsp
  &\ge\dd\frac12\,
  |\nabla\beta(t,u)|^2_2
  +\vp|\beta(t,u)|^2_2- C|b^*(t,u)|^2_2\vsp
  & \ge\dd\frac12\,|\nabla\beta(t,u)|^2_2+\vp|\beta(t,u)|^2_2- C|u|^2_2,\ \ff\,u\in D_0.\earr\hspace*{-5mm}
  \end{equation}
  Recall that we have the interpolation inequality
  $$|u|^2_2\le C|u|_{H^1}|u|_{-1},\ \ \ff\,u\in H^1.$$
  \n(Here, we have denoted by $C$ several positive constants depending on $\vp$, but independent of $u$.)}

 Noting that, by \eqref{e1.3} and \eqref{e1.4az},
 \begin{equation}\label{e3.12prim}
 \barr{lcl}
 |\nabla\beta(t,u)|^2_2&=&
 |\beta_x(t,u)+\beta_r(t,u)\nabla u|^2_2
 \ge\dd\frac12\,\nu^2|\nabla u|^2_2-|h|^2_\9|u|^2_2\\
 &\ge&\dd\frac14\,\nu^2|\nabla u|^2_2-C|u|^2_{-1},\earr\end{equation} 	
   \eqref{e3.9bbb} yields
  $$|u|^2_1+|A_\vp(t,u)|^2_{-1}\ge C(|\beta(t,u)|^2_2+|\nabla\beta(t,u)|^2_2
  +|b^*(t,u)|^2_2,\ \ff u\in D_0,$$and so, by \eqref{e3.9bb}, we get
  $$|A_\vp(t)u-A_\vp(s)u|_{-1}\le C|t-s|(|A_\vp(t)u|_{-1}+|u|_{-1}+1),\ \ff(t,s)\in[0,T],\,u\in D_0.$$This yields
   $$|A_\vp(t)u-A_\vp(s)u|_{-1}\le  C|t-s|L(|u|_{-1})(1+|A_\vp(t)u|_{-1}),\ \ff\,t,s\in[0,T],\,u\in D_0,$$
   where $L:[0,\9)\to[0,\9)$ is defined by
   $$L(r)=\max\left\{1,r+1\right\},\ \ \ff\,r\ge0.$$ Hence \eqref{e3.9a} follows.

\mk\noindent{\bf Proof of Proposition \ref{p2.4}.} By  Lemma \ref{l3.1} and \cite[Theorem 3.4]{5},    it follows that, for each $u_0\in D_0$,   there is a unique strong solution $u=u_\vp\in C([0,T];H^{-1})$ to equation \eqref{e2.7}, given by the exponential formula
\begin{equation}\label{e3.5a}
u_\vp(t)=\lim_{n\to\9}\prod^n_{k=1}\(I+\frac{t}n\,A_\vp\(k\,\frac{t}n\)\)^{-1}u_0,\ \ff t\in[0,T],
\end{equation}and \eqref{e3.5a} is uniform in $t$,
where the limit is taken in the strong topology of $H\1$. As a matter of fact, \eqref{e3.5a} is just the   finite difference scheme corres\-pon\-ding to the Cauchy problem \eqref{e2.7}. Namely,
\begin{equation}\label{e3.8a}
\barr{l}
\mu\1(u^\vp_{i+1}-u^\vp_{i})+A_\vp((i+1)\mu)u^\vp_{i+1}=0,\ i=0,1,...,N=\mbox{$\left[\frac T\mu\right]$}\vsp
u^\vp_0=u_0\earr\end{equation}
\begin{equation}\label{e3.8aa}
\dd\lim_{\mu\to0}u^\vp_\mu(t)=u_\vp(t)\mbox{ strongly in }H^{-1},\ \ff t\in[0,T],\end{equation}where for   $\mu>0$  the step function $u^\vp_\mu:[0,T]\to H\1$ is defined by
\begin{equation}\label{e3.8aaa}
u^\vp_\mu(t)=u^\vp_i,\ \ff t\in[i\mu,(i+1)\mu). \end{equation}
  Moreover, since $u_0\in D_0$, then   we have  $u_\vp(t)\in D_0$ for every $t\in[0,T]$ and
\begin{eqnarray}
&\dd\frac{du_\vp}{dt}\in L^\9(0,T;H\1),\ A_\vp u_\vp \in L^\9(0,T;H\1),\label{e3.18a1}\\
&\dd u_\vp(t),\beta(t,u_\vp(t))\in H^1,\ \ \ff t\in[0,T].\ \ \label{e3.18a}
\end{eqnarray}
%where \eqref{e3.18a1} follows from \eqref{e3.9a} by similar arguments as in the proof of \mbox{\cite[Theorem 4.19, p.~182]{1}}.

By \eqref{e3.18a1}--\eqref{e3.18a}, (H2) and  \eqref{e2.6}, it follows that
$$\barr{l}
\dd\frac12\ \frac d{dt}\,|u_\vp(t)|^2_{-1}+
\left<\beta(u_\vp(t)),u_\vp(t)\right>_2 -\left<{\rm div}(b(t, u_\vp)u_\vp(t),u_\vp(t))\right>_{-1}\vsp
\qquad\le|u_\vp(t)|_{-1}
|{\rm div}(b(t, u_\vp(t))u_\vp(t))|_{-1}\le C|u_\vp(t)|_{-1}
|b(t, u_\vp(t))u_\vp(t)|_2\vsp
\qquad\le C_1|u_\vp(t)|_{-1}|u_\vp(t)|_2,\mbox{\ \ a.e. }t\in(0,T).\earr$$By \eqref{e1.3}, this yields
$$\frac d{dt}\,|u_\vp(t)|^2_{-1}+\nu|u_\vp(t)|^2_2\le C_2|u_\vp(t)|^2_{-1},\mbox{\ \ a.e. }t\in(0,T),$$and, therefore, by Gronwall's lemma
\begin{equation}
\label{e3.20d}
|u_\vp(t)|^2_{-1}+\int^T_0|u_\vp(t)|^2_2dt\le C_3|u_0|^2_{-1},\ \ff\,t\in[0,T].\end{equation}Next, we have
$$\barr{r}
\beta(u_\vp(t))=(\vp I-\Delta)\1(A_\vp(t)u_\vp(t))+(\vp I-\Delta)\1{\rm div}(b(t,u_\vp(t))u_\vp(t)),\vsp\mbox{ a.e. }t\in(0,T),\earr$$and since, by \eqref{e3.18a1}, $A_\vp u_\vp\in L^2(0,T;H\1)$, we have
$$(\vp I-\Delta)\1 A_\vp u_\vp\in L^2(0,T;H^1),$$while
$$\barr{r}
|(\vp I-\Delta)\1{\rm div}(b(t,u_\vp(t))u_\vp(t))|_{H^1}\le C|b(t,u_\vp(t))u_\vp(t)|_2\le C_1|u_\vp(t)|_2,\vsp\mbox{ a.e. }t\in(0,T).\earr$$ \n Then, by \eqref{e3.20d}, we infer that
\begin{equation}
\label{e3.20dd}
\beta(u_\vp),u_\vp\in L^2(0,T;H^1).
\end{equation}
\n This completes the proof of Proposition \ref{p2.4}.

\bigskip\noindent{\bf Proof of Theorem \ref{t2.1}.} In the following we shall omit $x$ in the notations $\beta(t,x,u)$ and $b(t,x,u)$. For $\vp>0$, we consider the solution $u_\vp$ to \eqref{e2.6} and get first a few apriori estimates. In this proof, we consider $H^{-1}$ with its usual inner product $\left<u,v\right>_{-1}=\left<I-\Delta)^{-1}u,v\right>_{L^2}$. We have

\begin{equation}\label{e3.7a}
\barr{l}
\dd|u_\vp(t)|^2_2+\int^t_0\int_{\rr^d}(|\nabla u_\vp(s,x)|^2+|\nabla\beta(s,x,u_\vp(s,x))|^2)ds\,dx\vsp\qquad
+\vp\dd\int^t_0\int_{\rr^d}|\beta(s,x,u_\vp)|^2ds\,dx
\le C_T|u_0|^2_2,\ \ff t\in[0,T].\earr\end{equation}

To prove \eqref{e3.7a}, we need the following:

For some $C$ $(=C(|a|_\9,|b|_\9,\nu,T)\in(0,\9)),$
\begin{equation}\label{e3.24r}
	|u_\vp(t)|^2_2+\int^t_0|\nabla u_\vp(s)|^2_2ds
	+\vp\int^t_0|\beta(s,u_\vp(s))|^2_2ds
	\le C|u_0|^2_2,\ t\in[0,T],
	\end{equation}which, by \eqref{e3.18a1} and \cite[Theorem 1.19, p.~25]{1} in turn implies
	\begin{equation}\label{e3.25r}
	u_\vp\in C([0,T];L^2).
	\end{equation} 
To prove \eqref{e3.24r}, we note that by \eqref{e3.18a1}, \eqref{e3.20dd} in \eqref{e2.6} (with $u$ replaced by~$u_\vp$) we can take $L^2$-inner product with $u_\vp$ and integrate over $(0,t)$ to obtain for every $t\in[0,T]$
$$\barr{c}
\dd\frac12\,|u_\vp(t)|^2_2
+\dd\int^t_0\<\nabla u_\vp(s),\nabla\beta(s,u_\vp(s))\>_2ds
+\vp\dd\int^t_0\<u_\vp(s),\beta(s,u_\vp(s))\>_2ds\\
 =\dd\frac12\,|u_0|^2_2+\dd\int^t_0\<\nabla u_\vp(s),b^*(s,u_\vp(s))\>_2ds,\earr$$where the first integral on the left hand side by \eqref{e1.3}, \eqref{e1.4az} is bigger than
$$\int^t_0\(\frac{3\nu}4\,|\nabla u_\vp(s)|^2_2-\frac1\nu\,|h|^2_\9|u_\vp(s)|^2_2\)ds$$and the integral on the right hand side by (H2) is dominated by
$$\int^t_0\(\frac{\nu}4\,|\nabla u_\vp(s)|^2_2+\frac1\nu\,|b|^2_\9|u_\vp(s)|^2_2\)ds.$$Hence, since by (H1) $|\beta(s,u_\vp(s))|\le|a|_\9|u_\vp(s)|$, by Gronwall's inequality we obtain \eqref{e3.24r}.

Furthermore, we need that, for a.e. $t\in[0,T]$,
\begin{equation}\label{e3.26r}
\barr{ll}
\dd\frac d{dt}\int_{\rr^d}j(t,x,u_\vp(t,x))dx\!\!&={\raise-3mm\hbox{$_{H^1}$}}\!\!\<\beta(t,u_\vp(t)),\frac{du_\vp(t)}{dt}\>_{H^{-1}}\vsp &
+\dd\int_{\rr^d}j_t(t,x,u_\vp(t,x))dx,\earr\end{equation}
where $j(t,x,r)=\int^r_0\beta(t,x,\bar r)d\bar r,\ r\in\rr.$ 
To prove \eqref{e3.26r}, we first note that due to (H1), for all $(t,x,r)\in[0,T]\times\rr^d\times\rr$,
\begin{equation}
\label{e3.27r}
\nu|r|^2\le j(t,x,r)\le|a|_\9|r|^2,|j_t(t,x,r)|\le\frac12\,h(x)|r|^2.
\end{equation}Fix $t\in[0,T]$ and define the convex lower semi-continuous function $\vf^t:H\1\to[0,\9]$ by
	$$\vf^t(u):=\left\{\barr{ll}
	\dd\int_\rrd j(t,x,u)dx&\mbox{ if }u\in L^2,\vsp
	+\9&\mbox{ if }u\in H\1\setminus L^2.\earr\right.$$
	Then it is elementary to check that for its subdifferential $\pp\vf^t$ on $H\1$ we have
	$$(I-\Delta)\beta(t,u)=\pp\vf^t(u),\ \ff u\in L^2.$$
	Furthermore, \eqref{e3.20dd} and \eqref{e1.3prim}--\eqref{e1.4az} imply that
	$$\beta(t,u_\vp)\in L^2(0,T;H^1).$$Hence we conclude by \cite[Lemma 4.4, p.~158]{1} that
	\begin{equation}
	\label{e3.28r}
	\frac d{ds}\,\vf^t(u_\vp(s))=
{\raise-3mm\hbox{$_{H^1}$}}\!\!\<\beta(t,u_\vp(t)),
\frac{du_\vp(s)}{ds}\>_{H^{-1}}	\mbox{ for }ds\mbox{-a.e. }s\in[0,T].	
	\end{equation}
	Applying \eqref{e3.28r} with $s=t$, we find for $dt$-a.e. $t\in[0,T]$
		\begin{equation}
	\label{e3.29r}
	\barr{l}
	\dd\frac d{dt}\int_\rrd j(t,x,u_\vp(t,x))dx\\
	 \qquad=\dd\lim_{\delta\to0}\Big[\frac1\delta\int_\rrd(j(t+\delta,x,u_\vp(t+\delta,x))
	-j(t,x,u_\vp(t+\delta,x)))dx\\
	\qquad\qquad\qquad\qquad+\ \dd\frac1\delta(\vf^t(u_\vp(t+\delta))-\vf^t(u_\vp(t)))\Big]\vsp
	 \qquad=\dd\lim_{\delta\to0}\int_\rrd\int^1_0j_t(t+s\delta,x,u_\vp(t+\delta,x))ds\,dx\\\qquad\qquad\qquad\qquad
	+{\raise-3mm\hbox{$_{H^1}$}}\!\!\<\beta(t,u_\vp(t)),
	\dd\frac{du_\vp(t)}{dt}\>_{H^{-1}}.\earr	
	\end{equation}
	Since, as $\delta\to0$,
	$$\mbox{$j_t(t+s\delta,x,u_\vp(t+\delta,x))\to
		j_t(t,x,u_\vp(t,x))$ in $dt\otimes dx$ measure,}$$ and by \eqref{e3.27r}, \eqref{e1.4az}
	 $$j_t(t+s\delta,x,u_\vp(t+\delta,x))\le\frac12\,h(x)|u_\vp(t+\delta,x)|^2$$for all $\delta<\delta_0$ and all $s\in[0,1]$, $(t,x)\in[0,T]\times\rrd$, where the latter term by \eqref{e3.25r} converges in $L^1(\rrd)$ as $\delta\to0$, \eqref{e3.29r} implies \eqref{e3.26r}.
	
	Now, it is easy to prove \eqref{e3.7a}. Integrating \eqref{e3.26r} over $(0,t)$, by \eqref{e2.6} and \eqref{e3.27r}, we find		
	$$\barr{lcl}
	0&\le&\dd\int_\rrd j(t,x,u_\vp(t,x))dx\vsp&\le&\dd\int_\rrd j\(0,x,u_0(x)\)dx+\frac12\,|h|_\9\int^t_0|u_\vp(s)|^2_2ds\vsp&&-\dd\int^t_0|\nabla\beta(s,u_\vp(s))|^2_2-\vp\int^t_0|\beta(s,u_\vp(s))|^2_2ds\vsp&&+\dd\int^t_0\<b^*(s,u_\vp(s)),\nabla\beta(s,u_\vp(s))\>_2dx,\earr$$which by (H1), (H2) and \eqref{e3.24r}, \eqref{e3.27r} implies
	 $$0\le\frac12\,(|a|_\9+(|h|_\9+|b|^2_\9)T)|u_0|^2_2-\frac12\int^t_0|\nabla\beta(s,u_\vp(s))|^2_2ds.$$Therefore, by \eqref{e3.24r} inequality \eqref{e3.7a} follows. We note also that, by \eqref{e2.6} and \eqref{e3.7a}, it follows that  $\left\{\frac{du_\vp}{dt}\right\}$ is bounded in $L^2(0,T;H\1)$.
Hence, along a subsequence, again denoted $\{\vp\}\to0$, we have
$$\barr{rcll}
u_\vp&\longrightarrow&u&\mbox{weak-star in $L^\9(0,T;L^2)$}\\
&&&\mbox{weakly in $L^2(0,T;H^1)$}\vsp
\beta(u_\vp)&\longrightarrow&\eta&\mbox{weakly in $L^2(0,T;H^1)$}\vsp
\dd\frac{du_\vp}{dt}&\longrightarrow&\dd\frac{du}{dt}&\mbox{weakly in $L^2(0,T;H\1)$}\vsp
b(u_\vp)u_\vp&\longrightarrow&\zeta&\mbox{weakly in $L^2(0,T;L^2)$}\earr$$where
 \begin{equation}\label{e3.8} %{e3.30}%
 \barr{c}
 \dd\frac{du}{dt}-\Delta\eta+{\rm div}\ \zeta=0\mbox{ in }\cald'((0,T)\times\rr^d),\\
 u(0,x)=u_0(x),\ \ x\in\rr^d.\earr\end{equation}
 Moreover, since $H^1$ is compactly embedded in $L^2_{\rm loc}(\rr^d)$, by the Aubin-Lions compactness theorem (see, e.g., \cite{1}, p.~26), it follows that, for $\vp\to0$, we have
 \begin{equation}\label{e3.30b}
 	\barr{rcll}
 u_\vp&\longrightarrow&u&\mbox{strongly in $L^2(0,T;L^2_{\rm loc}(\rr^d))$},\earr\end{equation}hence selecting another subsequence, if necessary, $u_\vp\to u$, a.e., and so
 $$\beta(u_\vp) \longrightarrow \beta(u),\ \mbox{a.e.}$$
 Hence $\eta=\beta(u)$ and $\zeta=b(u)u$, a.e. on $(0,T)\times\rr^d$ and thus $u$ solves \eqref{e1.1}

 Letting $\vp\to0$ in \eqref{e3.7a} and taking into account \eqref{e2.6}, we obtain the estimates
 $$\barr{c}
 \dd\|u(t)\|^2_{L^\9(0,T;L^2)}+\int^T_0\int_{\rr^d}
 (|\nabla u(t,x)|^2+|\nabla\beta(t,x,u(t,x))|^2)dt\,dx\le C|u_0|^2_2,\vsp
 \left\|\dd\frac{du}{dt}\right\|_{L^2(0,T;H^{-1})}\le C.\earr$$

 Hence $u\in L^\9(0,T;L^2)\cap L^2(0,T;H^1)\cap W^{1,2}([0,T];H\1)$ and $\beta(u)\in  L^2(0,T;H^1)$, in particular, by \cite[Theorem 1.19, p.~25]{1} $u\in C([0,T];L^2)$.

 This implies that $u$ satisfies \eqref{e2.1}, \eqref{e2.2} and so it  is a strong solution to \eqref{e1.1}. The uniqueness of the solutions $u$ satisfying \eqref{e2.1}--\eqref{e2.2} is based on the monotonicity of \mbox{$r\to\beta(\cdot,\cdot,r)$} and the fact that   $b^*$ is Lipschitz by (H2). Indeed, if $u_1,u_2$ are two strong solutions to \eqref{e1.1}, we have
 $$\barr{l}
 \dd\frac12\ \frac d{dt} \ |u_1(t)-u_2(t)|^2_{-1}+
 (\beta(t,u_1(t))-\beta(t,u_2(t)),u_1(t)-u_2(t))_2\vsp
\hspace*{23mm} +\<{\rm div}(b(t,u_1)u_1-b(t,u_2)u_2),u_1-u_2\>_{-1}\vsp
 \hspace*{23mm} =\<\beta(t,u_1)-\beta(t,u_2),u_1-u_2\>_{-1}.\earr$$
 Then, by (H1)--(H2), we get
 $$\barr{l}
 |u_1(t)-u_2(t)|^2_{-1}+\dd\int^t_0|u_1(s)-u_2(s)|^2_2\vsp
 \quad\le C_1\!\dd\int^t_0(|b^*(s,u_1){-}b^*(s,u_2)|^2_2
 {+}|\beta(s,u_1){-}\beta(s,u_2)|_2
 |u_1(s){-}u_2(s)|_{-1})ds\vsp
 \quad\le C_2\!\dd\int^t_0(|u_1(s)-u_2(s)|^2_2+|u_1(s)-u_2(s)|^2_{-1})ds,\ \ff t\in(0,T),
 \earr$$and so $u_1\equiv u_2$, as claimed. \newpage 

 Assume now that $u_0\in L^1\cap D_0.$ To prove that $u$ satisfies \eqref{e1.5}--\eqref{e1.7}, we multiply equation \eqref{e1.1} by  $\calx_\delta(u(t,x))$, where $\delta>0$ and $\calx_\delta$ is the following approximation of the signum function
 $$\calx_\delta(r)=\left\{\barr{cl}
 1&\mbox{ for }r\ge\delta,\\
 \dd\frac r\delta&\mbox{ for }|r|<\delta,\\
 -1&\mbox{ for }r\le-\delta.\earr\right.$$
 We note that, for $dt$-a.e. $t\in(0,T)$, since $u(t)\in H^1(\rr^d)$,  it follows that $\calx_\delta(u(t))\in H^1(\rr^d)$.
If we apply ${\raise-1mm\hbox{$_{H^{-1}}$}}\!\<\cdot,\calx_\delta(u(t,\cdot))\>_{H^1}$ to \eqref{e1.1}   and integrate over $(0,t)$, we get
 $$\barr{l}
\dd\int_{\rr^d}j_\delta(u(t,x))dx
+\int^t_0\int_{\rr^d}\nabla\beta(s,x,u(s,x))
\cdot\nabla\calx_\delta(u(s,x))dx\,ds\vsp
\qquad=\dd\int^t_0\!\!\int_{\rr^d}
b^*(s,u(s,x))
\cdot\nabla\calx_\delta(u(s,x))dx\,ds 
+\dd\int_{\rr^d}j_\delta(u_0(x))dx,\ \ff t\in[0,T],\earr$$where
\begin{equation}\label{e3.31b}
j_\delta(r)=\int^r_0\calx_\delta(s)ds,\ \ff r\in\rr.\end{equation}
Taking into account that by \eqref{e1.3} $\beta_u,$ $\calx'_\delta\ge0,$ a.e. on $(0,T)\times\rr^d$ and $\rr$, respectively, we have  by \eqref{e1.4az}  that, for $dt$-a.e. $t\in(0,T)$,
\begin{equation}\label{e3.10a}
	\hspace*{-3mm}\barr{ll}
\dd\int_{\rr^d}\!\!\nabla\beta(t,x,u(t,x))\cdot\nabla\calx_\delta(u(t,x))dx\vsp
\qquad
= \dd\int_{\rr^d}\!(\beta_u(t,x,u(t,x))
\nabla u(t,x)+\beta_x(t,x,u(t,x)))
\nabla \calx_\delta(u(t,x))dx\vsp
\qquad\ge\dd\int_{\rr^d}\b_x(t,x,u(t,x))\cdot\nabla\calx_\delta(u(t,x))dx
\vsp\qquad\dd\ge-\frac 1\delta\int_{[|u(t,x)|\le\delta]}h(x)|u(t,x)|\cdot|\nabla u(t,x)|dx\vsp\qquad
\ge-|h|_2\(\dd\int_{[|u(t,x)|\le\delta]}|\nabla u(t,x)|^2dx\)^{\frac12}
=-\eta_\delta(t).\earr\hspace*{-3mm}\end{equation}
Since $\nabla u(t,x)=0$, a.e. on $[x\in\rr^d;\ u(t,x)=0]$, it follows that, for  $dt$-a.e. $t\in(0,T)$,   $\eta_\delta(t)\to0$ as $\delta\to0$.
We have, therefore,
$$\liminf_{\delta\to0}\int^t_0\int_{\rr^d}\nabla\beta(s,x,u(s,x))\cdot\nabla\calx_\delta(u(s,x))dx\,ds\ge0,$$and, similarly, it follows by (H2) part \eqref{e1.7prim} that
\begin{equation}
	\label{e3.10b}
	\barr{r}
	\dd\lim_{\delta\to0}\left|
 \dd\int^t_0\int_{\rr^d}b^*(s,u(s,x))\cdot\nabla\calx_\delta(u(s,x))\right|dx\,ds\vsp \quad\le\dd\lim_{\delta\to0}\int^t_0\int_{[|u(s,x)|\le\delta]}h(x)|\nabla u(s,x)|dx\,ds=0,\earr\end{equation}because $h\in L^2$. This yields, since $0\le j_\delta(r)\le|r|$, $r\in\rr$,
$$ \liminf_{\delta\to0} \int_{\rr^d}j_\delta(u(t,x))dx
\le\int_{\rr^d}|u_0(x)|dx$$and, since $\lim\limits_{\delta\to0}j_\delta(r)=|r|$,   $r\in\rr$, we infer by Fatou's lemma that
 \begin{equation}\label{e3.9}
 \int_{\rr^d}|u(t,x)|dx\le\int_{\rr^d}|u_0(x)|dx,\ \ff t\in[0,T],\end{equation}and, therefore,  $u\in L^\9(0,T;L^1(\rr^d))$ with $|u(t)|_1\le|u_0|_1$, for all $t\in[0,T]$.

 If $u_0\in D_0$ and $u_0\ge0$, a.e. on $\rr^d$, then applying
 ${\raise-1mm\hbox{$_{H^{-1}}$}}\!\<\cdot,u^-(t)\>_{H^1}$ to  \eqref{e1.1}  and integrating over $(0,T)$, we see that

 $$\int_{\rr^d}|u^-(t,x)|^2dx=0,\ \ff t\in[0,T],$$and so $u\ge0$, a.e. in $(0,T)\times\rr^d$ and \eqref{e1.5} is proved.

 Finally, integrating \eqref{e1.1} over
 $\rr^d$, we see that \eqref{e1.7} holds. (More precisely, one multiplies \eqref{e1.1} with $\psi_n\in C^2_0(\rr^d)$, $n\in\nn$, with uniformly bounded first and second order derivatives such that $0\le\psi_n\nearrow 1$, integrates over $\rr^d$ and lets $n\to\9.$)

 It remains to prove the weak-continuity condition \eqref{e1.6}. So, let $u_0\in L^1\cap D_0$ with $u_0\ge0$ and $|u_0|_1=1.$ Let $u$ be the solution to \eqref{e1.1} with the initial condition $u_0$ and let $t_n$, $t\in[0,T]$, such that $t_n\to t$ as  $n\to\9$. Then, the probability measures $\mu_n(dx):=u(t_n,x)dx$ converge vaguely to the probability measure $\mu(dx):=u(t,x)dx$, because $u\in C([0,T];L^2)$. Since, for probability measures, vague and weak convergence are equivalent, \eqref{e1.6} follows.

 Assume now \eqref{e2.3prim}   and let us prove \eqref{e2.3a}. To this end, it is convenient to use the finite difference scheme \eqref{e3.8a}-\eqref{e3.8aa}. Namely,
 
\begin{equation}\label{e3.35b}
\barr{l}
\dd\frac1\mu\,(u^\vp_{i+1}-u^\vp_i)-\Delta\beta((i+1)\mu,u^\vp_{i+1}){+}\vp\beta((i+1)\mu,u^\vp_{i+1})\\\qquad\quad+{\rm div}(b^*((i+1)\mu,u^\vp_{i+1}))=0,\  i=0,1,...,N\mbox{ in }\rr^d,\vsp
u^\vp_0=u_0,\earr\end{equation}
for the solution $u^\vp=u^\vp(t,u_0)$ to \eqref{e2.6} and, similarly, for $\bar u^\vp=u^\vp(t,\bar u_0)$, where $u_0,\bar u_0\in L^1\cap D_0$. Here again we supress the $x$-dependence in the notation. We have
\begin{equation}\label{e3.34b}\int_{\rr^d}|u^\vp_{i+1}|dx\le\int_{\rr^d}|u_0|dx,\ \ \ff\,i=0,1,...,N.\end{equation}We postpone for the time being the proof of \eqref{e3.34b} and
  we set $u_{i+1}=u^\vp_{i+1}$ and $\bar u_{i+1}=\bar u^\vp_{i+1}$ (corresponding to $\bar u_0$).  Also we shall write
  $$\beta_i(u_i)=\beta(i\mu,x,u_i(x)),\ \beta_i(\bar u_i)=\beta(i\mu,x,\bar u_i(x)),$$and recall also that $b^*(t,x,u)$ is simply denoted $b^*(t,u)$. We get
\begin{equation}\label{e3.34bb}
\hspace*{-4mm}\barr{l}
\dd\frac1\mu(u_{i+1}{-}\bar u_{i+1}){-}\Delta(\beta_{i+1}(u_{i+1}){-}\beta_{i+1}(\bar u_{i+1})){+}\vp(\beta_{i+1}(u_{i+1}){-}\beta_{i+1}(\bar u_{i+1}))\vsp
\qquad\qquad +\,{\rm div}(b^*(({i+1})\mu,u_{i+1}){-}b^*(({i+1})\mu,\bar u_{i+1}))=\dd\frac1\mu(u_i{-}\bar u_i),\vsp\hfill i=0,1,...,N.
\earr\end{equation}
Then, multiplying by $\calx_\delta(\beta_{i+1}(u_{i+1})-\beta_{i+1}(\bar u_{i+1}))$ and integrating over $\rr^d$, we~get
\begin{equation}\label{e3.33a}
\hspace*{-6mm}\barr{l}
\dd\frac1\mu\int_{\rr^d}(u_{i+1}{-}\bar u_{i+1})\calx_\delta(\beta_{i+1}(u_{i+1}){-}\beta_{i+1}(\bar u_{i+1}))dx\vsp
\qquad+\vp
\dd\int_{\rr^d}(\beta_{i+1}(u_{i+1})-\beta_{i+1}(\bar u_{i+1}))
\calx_\delta(\beta_{i+1}(u_{i+1}){-}\beta_{i+1}(\bar u_{i+1}))dx\vsp
\qquad+\dd\int_{\rr^d}|\nabla(\beta_{i+1}(u_{i+1})-\beta_{i+1}(\bar u_{i+1}))|^2\calx'_\delta(\beta_{i+1}(u_{i+1})-\beta_{i+1}(\bar u_{i+1}))dx\vsp
\qquad-\dd\int_{\rr^d}(b^*(({i+1})\mu,u_{i+1}){-}b^*(({i+1})\mu,\bar u_{i+1}))\cdot\vsp\hfill\cdot\nabla(\beta_{i+1}(u_{i+1}){-}\beta_{i+1}(\bar u_{i+1}))\calx'_\delta(\beta_{i+1}(u_{i+1}){-}\beta_{i+1}(\bar u_{i+1}))dx\vsp\qquad=\dd\frac1\mu\int_{\rr^d}(u_i-\bar u_i)\calx_\delta(\beta_{i+1}(u_{i+1})-\beta_{i+1}(\bar u_{i+1}))dx.
\earr\hspace*{-6mm}\end{equation}
We set
$$
v^\delta_{i+1}=\calx_\delta(\beta_{i+1}(u_{i+1})-\beta_{i+1}(\bar u_{i+1}))-\calx_\delta(u_{i+1}-\bar u_{i+1})$$
and note that, clearly,
\begin{equation}
\label{e3.33aa}
\sup_{\delta\in(0,1)}\|v^\delta_{i+1}\|_\9\le2\mbox{\ \  and\ \ }\lim_{\delta\to0}v^\delta_{i+1}=0,\mbox{ a.e. on }\rr^d,\end{equation}and
$$\barr{l}
((u_{i+1}-\bar u_{i+1})-(u_i-\bar u_i))\calx_\delta(u_{i+1}-\bar u_{i+1})\vsp\hfill
\ge j_\delta(u_{i+1}-\bar u_{i+1})-j_\delta(u_i-\bar u_i),\ \ff i,\earr$$where
$j_\delta$ is defined by \eqref{e3.31b}. Then, by \eqref{e3.33a}, this yields
\begin{equation}\label{e3.33aaa}
\barr{l}
\dd\int_{\rr^d}j_\delta(u_{i+1}-\bar u_{i+1})
dx\le
\dd\int_{\rr^d}j_\delta (u_0-\bar u_0) dx\\
\quad-\dd\sum^i_{j=0}\int_{\rr^d}v^\delta_{j+1}
[(u_{j+1}-\bar u_{j+1})-(u_j-\bar u_j)]dx\vsp
\quad+
\dd\sum^i_{j=0}\mu\int_{\rr^d}
(b^*((j+1)\mu,u_{j+1}){-}b^*((j+1)\mu,\bar u_{j+1}))\cdot\vsp\qquad\cdot\nabla
(\beta_{j+1}(u_{j+1}){-}\beta_{j+1}(\bar u_{j+1}))\calx'_\delta(\beta_{j+1}(u_{j+1}){-}\beta_{j+1}(\bar u_{j+1}))dx,
\earr\end{equation}for all $i\le N$.
On the other hand, we have by \eqref{e1.3} and \eqref{e2.3prim}
$$\barr{l}
\!\!\dd\lim_{\delta\to0}\int_{\rr^d}\!\!
|b^*((j+1)\mu,u_{j+1}){-}b^*((j+1)\mu,\bar u_{j+1})|\vsp\qquad\qquad
 \cdot
|\nabla(\beta_{j+1}(u_{j+1}){-}\beta_{j+1}(\bar u_{j+1}))|
\calx'_\delta(\beta_{j+1}(u_{j+1}){-}\beta_{j+1}(\bar u_{j+1}))dx\vsp
\le \dd\lim_{\delta\to0}\,\dd\frac1\delta
\int_{\left[|u_{j+1}-\bar u_{j+1}|\le\frac\delta\nu\right]}h(x)
|u_{j+1}-\bar u_{j+1}|\vsp\qquad\qquad \cdot|\nabla(\beta_{j+1}(u_{j+1})-\beta_{j+1}(\bar u_{j+1}))|dx=0,
\earr$$
because $h,\nabla\beta_{j+1}(u_{j+1}),\nabla\beta_{j+1}(\bar u_{j+1})\in L^2$ and   \mbox{$\nabla(\beta_{j+1}(u_{j+1}){-}\beta_{j+1}(\bar u_{j+1}))=0$} on $[|u_{j+1}-\bar u_{j+1}|=0]$. Moreover, by
\eqref{e3.34b}, \eqref{e3.33aa} and Lebesgue's dominated convergence theorem,  we have
$$\lim_{\delta\to0}\sum^i_{j=0}\int_{\rr^d} v^\delta_{j+1}
[(u_{j+1}-\bar u_{j+1})-(u_j-\bar u_j)]dx=0,\ \ff\,i=1,...,N.$$Then, by \eqref{e3.34b} and \eqref{e3.33aaa}, we get
$$|u_{i+1}-\bar u_{i+1}|_1\le|u_0-\bar u_0|_1,\ \ff i=0,1,...,N.$$Then, by  \eqref{e3.8aaa}, we get
$$|u^\vp_h(t)-\bar u^\vp_h(t)|_1\le|u_0-\bar u_0|_1,\ \ff u_0,\bar u_0\in L^1\cap L^2,$$and so, by \eqref{e3.8aa}, it follows that
\begin{equation}\label{e3.13a}
|u_\vp(t,u_0)-u_\vp(t,\bar u_0)|_1\le|u_0-\bar u_0|_1,\ \ff u_0,\bar u_0\in L^1\cap L^2,\ \ff t\in(0,T].\end{equation}

  Moreover, by \eqref{e3.30b} it follows that, for $\vp\to0$, $u_\vp(\cdot,u_0)\to u(\cdot,u_0)$ in $L^2(0,T;B_R)$ for all $B_R=\{x\in\rr^d;\|x|\le R\}$ and, therefore, \eqref{e3.13a} yields
$$\|u(t,u_0)-u(t,\bar u_0)\|_{L^1(B_R)}\le|u_0-\bar u_0|_1,\ \mbox{ a.e. }t\in(0,T),$$
for all $R>0$ and so \eqref{e2.3a} follows.

\bk\n{\bf Proof of \eqref{e3.34b}.} Taking into account that
$$\calx_\delta(u_{i+1})(u_{i+1}-u_i)\ge j_\delta(u_{i+1})-j_\delta(u_i),\ \ff\,i=0,1,...,N,\mbox{ a.e. in }\rr^d,$$and arguing as in \eqref{e3.33a}--\eqref{e3.33aaa}, we see that
$$\barr{l}
\dd\int_{\rr^d}j_\delta(u_{i+1})dx\le\dd\int_{\rr^d}|u_0|dx+\mu\sum^i_{j=0}\int_{\rr^d}b^*((j+1)\mu,u_{j+1})\cdot\nabla u_{j+1}\calx'_\delta(u_{j+1})dx\vspace*{-2mm}\\
\qquad\hfill \ff\,i=0,1,...,N,\earr$$which, as above, by \eqref{e1.7prim} for $\delta\to0$ yields \eqref{e3.34b}, as claimed.

Assume now that   $u_0\in L^1\cap D_0\cap L^\9$ and that \eqref{e2.5prim} holds. If $u_\vp$ is the solution to \eqref{e2.7}, then we have by \eqref{e3.8a}--\eqref{e3.8aaa} that 
\begin{equation}\label{e3.43}
u_\vp(t)=\lim\limits_{\mu\to0}u^\mu_\vp(t)\ \mbox{ in $H\1$,}\ \ff\,t\in[0,T],\end{equation} where $u^\mu_\vp(t)=u_i$, $\ff\,t\in[i\mu,(i+1)\mu)$ and $\{u_i\}$ are given by
$$\mu\1(u_{i+1}-u_i)+(\vp I-\Delta)\beta_{i+1}(u_{i+1})+{\rm div}\,b^*((i+1)\mu,u_{i+1})=0,\ i=0,1,...,N.$$
Then, by \eqref{e3.34b}, we have that $u_i\in L^1$, $\ff i$. Then, by \eqref{e2.5prim}, it follows that, for $M=\Lambda(b,\beta),$ we have, since $\beta_{i+1}(|u_0|_\9+M(i+1)\mu)\in H^1$ by \eqref{e1.4az} and $b^*((i+1)\mu,|u_0|_\9+M(i+1)\mu)\in L^2$ by \eqref{e1.7prim},\newpage
$$\barr{l}
u_{i+1}-u_i-M\mu+\mu(\vp I-\Delta)(\beta_{i+1}(u_{i+1})-\beta_{i+1}
(|u_0|_\9+M(i+1)\mu))\vsp\qquad+\mu\,{\rm div}(b^*((i+1)\mu,u_{i+1})-b^*((i+1)\mu,|u_0|_\9+M(i+1)\mu))\le 0\mbox{ in }\rr^d,\vsp\hfill i=0,1,...,N.
\earr$$
Multiplying the latter by $\calx_\delta(\beta_{i+1}(u_{i+1})-\beta_{i+1}(|u_0|_\9+M(i+1)\mu))^+)$, summing over $i$, integrating over $\rr^d$ and using  
\eqref{e3.33aa},  we get as above
$$\barr{r}\dd\ov{\lim_{\delta\to0}}\int_{\rr^d}
(u_{i+1}-u_0-M(i+1)\mu)^+\calx_\delta((u_{i+1}-|u_0|_\9-M(i+1)\mu))^+dx\le0,\\ \ff i=0,1,...,N.\earr$$
Hence, $$u_{i+1}\le|u_0|_\9+M(i+1)\mu,\mbox{ a.e. on }\rr^d,\ i=0,1,...,N,$$and, therefore,
$$u^\mu_\vp(t,x)\le MT+|u_0|_\9,\mbox{ a.e. in }(t,x)\in(0,T)\times\rr^d.$$Similarly, it follows that $$u^\mu_\vp(t,x)\ge-MT-|u_0|_\9,\mbox{ a.e. in }(t,x)\in(0,T)\times\rr^d.$$ We have, therefore, $$\|u^\mu_\vp(t,x)\|_{L^\9((0,T)\times\rr^d)}\le MT+|u_0|_\9,\ \ff\,\vp>0,\ \mu>0,$$and, taking into account \eqref{e3.43}, we infer that $u_\vp\in L^\9((0,T)\times\rr^d)$ and
\begin{equation}\label{e3.44}
	\|u_\vp\|_{L^\9((0,T)\times\rr^d)}\le MT+|u_0|_\9,\ \ff\vp>0.\end{equation}
	
Now, recalling \eqref{e3.30b}, it follows by \eqref{e3.44}   that $u\in L^\9((0,T)\times\rr^d)$ and $\|u\|_{L^\9((0,T)\times\rr^d)}\le MT+|u_0|_\9$, as claimed.~$\Box$ 

\begin{remark}\label{r3.2}\rm  Inequality \eqref{e3.2} is absolutely necessary for  the existence  of a strong solution $u$ in the space $H\1(\rr^d)$  to the non\-linear evolu\-tion equation \eqref{e2.7} and in our case and, as seen above, this is \mbox{implied} by (H1), (H2). However, if $A(t)$ is monotone (which is not, however, the~case here), demicontinuous and coercive from $V$ to $V'$, where $V\subset H\subset V'$ and $V$ is a re\-flexive Banach space, $H$ is a Hilbert space and $V'$ is the dual of~$V$,\break the exis\-tence of a solution $u\in L^2(0,T;V)$ with $\frac{du}{dt}\in L^2(0,T;V')$ for\break $\frac{du}{dt}+A(t)u=0$ follows assuming only that $t\to A(t)u$ is measurable. (See,~e.g., \cite{1}, p.~177.) This result can be applied in our case to the appro\-xi\-ma\-ting  equation \eqref{e2.6} and so Hypotheses \eqref{e1.4az} and \eqref{e1.4a} are no longer ne\-ces\-sary with $H=H\1(\rr^d)$, \mbox{$V=L^2(\rr^d)$,} $V'=H^{-1}(\rr^d)$, to get a solution $u_\vp\in L^2(0,T;L^2(0,T;L^2(\rr^d)))\cap C([0,T];H^{-1}(\rr^d))$ to \eqref{e2.6} with   $\frac{du_\vp}{dt}\in $ $  L^2(0,T;H^{-2}(\rr^d))$.  However, under such weak assumptions, one \mbox{cannot} pass to the limit $\vp\to0$ in \eqref{e2.6} to~get~\eqref{e1.1}. 	 \end{remark}

\begin{remark}\label{r3.3}\rm It should be emphasized that, by taking into account the phy\-si\-cal significance, $L^1(\rr^d)$ is the natural space for the \FP\ equation \eqref{e1.1} and, as seen above, it follows by \eqref{e2.3a} that, for $u_0\in L^1(\rr^d)$, there is a distributional solution $u=u(t,u_0)$ to this equation (i.e. $u$ satisfies \eqref{e2.4}) which is weakly continuous in $t$ in the sense of \eqref{e1.6}. In fact, $u(t,u_0)=\lim\limits_{n\to\9}u(t,u^n_0)$ in $L^1$, $u^n_0\in L^1\cap L^2,$ $u^n_0\to u_0$ in $L^1$ as $n\to\9$. However, the uniqueness of such a solution remains open in this  case.
 In fact, the uniqueness proof in Theorem \ref{t2.1} is based on the regularity properties \eqref{e2.1}--\eqref{e2.2}, but the argument is not applicable to solutions with low regularity as is the case for distributional solutions. However, if $\beta(t,r)\equiv\beta(r)$, then the arguments of \cite{5b}, \cite{7az}, \cite{5bb} can be applied to derive the uniqueness of the distributional solution to \eqref{e1.1} which are $L^1$-valued narrowly  continuous in~$t$.\end{remark}

\begin{remark}
	\label{r3.4}\rm One might suspect that Theorem \ref{t2.1} can be extended to more general nonlinear \FP\ equations
	\begin{equation}
		\label{e3.4a}
		\barr{c}
		u_t-\dd\sum^d_{i,j=1}D^2_{ij}(a_{ij}(t,x,u)u)+{\rm div}(b(t,x,u)u)=0,\\
		u(0,x)=u_0(x),\ \ x\in\rrd,\earr\end{equation}
	where the matrix  $((a_{ij})_r(t,x,r)r+a_{ij}(t,x,r))^d_{i,j=1}$ is positive definite and $a_{ij},$ $b$ satisfy regularity conditions as in (H1)--(H2). In \cite{3a}, the authors proved in the autonomous case the existence of a mild (generalized) solution in $L^1$, via nonlinear semigroup theory, but  the arguments used in the proof of Theorem \ref{t2.1}, based on the fact that $I-\Delta$ is an isometry from $H^1$ to $H\1$, are not applicable for \eqref{e3.4a}. Of course, it remains the possibility mentioned in Remark \ref{r3.2}, that is, to represent \eqref{e3.4a} as a Cauchy problem in the variational pair $(H^1,H\1)$. 
	\end{remark}

\begin{remark}
	\label{r3.5}\rm One might ask if Theorem \ref{t2.1} can be extended to a measure initial data $u_0$. In the autonomous case, such a result was established in \cite{5aa} for $d=3$ and $\beta$ with polynomial growth, in particular for $\beta(r)\equiv ar^\alpha$, $\alpha\ge1$. One may suspect that this result remains true in the present case for $\beta\equiv\beta(r)$ and $b$ as in \cite{5aa}, but the problem is open.
\end{remark}

	Finally, let us prove the application   to the McKean--Vlasov SDE \eqref{e1.8}:

\bk\n{\bf Proof of Corollary \ref{c2.2}.}   Let $u_0\in L^1\cap D_0$, $u_0\ge0$ and $\int_{\rr^d}u_0dx=1$. Let $u$ be the strong solution from Theorem \ref{t2.1}. Then $\mu_t(dx):=u(t,x)dx,$ $t\in[0,T]$, are probability measures on $\rr^d$, weakly continuous in $t$, which solve \eqref{e2.4} and
$$\int^T_0\int_{\rr^d}(|a(t,x,u(t,x))|+|b(t,x,u(t,x))|)u(t,x)dx\,dt<\9,$$since both $a$ and $b$ are bounded. Now, the assertion follows by \cite[Section~2]{3a}. The case where we merely have $u_0\in L^1$ follows by the same arguments and Remark \ref{r2.1prim}.~$\Box$

\bk\n{\bf Acknowledgement.} This work was supported by the DFG through CRC 1283 and by a grant of the Ministry of Research, Innovation and Digitization, CNCS-UEFISCDI project  PN-III-P4-PCE-2021-0006 within PNCDI III.

\end{document}